\documentclass[10pt]{article}
\usepackage{url}
\usepackage{amsfonts,amsthm,amsmath}
\usepackage{amssymb}
\usepackage{authblk}
\usepackage{xcolor}
\usepackage{graphicx}
\usepackage[top=1.25in, bottom=1.5in, left=1.5in, right=1.5in]{geometry}

\newtheorem{Example}{Example}[section]

\title{An Analysis and Critique of the Scoring Method Used for Sport Climbing at the 2020 Tokyo Olympics}
\author{Michela J.\ Stinson\thanks{This research was supported by the Social Sciences and Humanities Research Council of Canada
(SSHRC) Joseph-Armand Bombardier Doctoral Fellowship}}
\affil{Department of Recreation and Leisure Studies\authorcr University of Waterloo\authorcr
Waterloo, Ontario, N2L 3G1, Canada
\authorcr
\texttt{mk.stinson@uwaterloo.ca}}
\author{Douglas R.\ Stinson\thanks{This research was supported by  
Natural Sciences and Engineering Research Council of Canada discovery grant RGPIN-03882.
}}
\affil{David R.\ Cheriton School of Computer Science\authorcr University of Waterloo\authorcr
Waterloo, Ontario, N2L 3G1, Canada
\authorcr
\texttt{dstinson@uwaterloo.ca}}

\date{\today}

\begin{document}
	
\maketitle

\begin{abstract}
Sport climbing was a new Olympic event introduced at the Tokyo 2020 Olympics. It was composed of three disciplines, and the final rankings were determined by computing
the product of each climber's rankings in the three disciplines, with the lowest score winning. In this paper, we compare this product-based scoring method with the more usual sum-based method. As well, we propose and analyze a new method based on taking the sum of the square roots of each climber's rankings. 
\end{abstract}

\section{Introduction: Sport Climbing and the 2020 Olympics}
\label{intro.sec}

While outdoor rock climbing has a diverse and lengthy global history, competitive sport climbing\footnote{In rock climbing, the term ``sport climbing'' has historically been employed to mean climbing on natural rock with permanently-placed protective gear. In this paper, we use ``sport climbing'' to refer to the multi-discipline event of competitive climbing done on artificial walls.} is relatively new, with the earliest records of competitions being held on artificial walls in the 1980s \cite{Kiewa}. The discipline of competitive sport climbing has since become increasingly bifurcated from more traditional rock climbing, not only though changes in ethics (i.e., concern with Leave No Trace and other moral constraints), purpose (i.e., immersion in and/or ``conquering'' of nature), and location (i.e., on natural vs.\ artificial walls), but also through the increase of active governing bodies, organizations, and institutional logics surrounding competitive sport climbing \cite{BR19,Kiewa}. Competitive sport climbing is seen as the ``rationalization'' or ``quantification'' of rock climbing \cite{Heywood,Kiewa}, an assertion that has been made uniquely visible by the choice to include sport climbing in the 2020 Tokyo Olympic Games, and through the subsequent choice of scoring format.

In 2015, the International Federation of Sport Climbing (IFSC) and the Japanese Mountaineering Association (JMA) proposed that sport climbing be added to the roster of the 2020 Tokyo Olympic Games \cite{Adie}. This proposal was then approved by International Olympic Committee (IOC) at the 129th IOC Session in Rio de Janeiro, allowing sport climbing into the 2020 Tokyo Olympic Games alongside skateboarding, baseball/softball, karate, and surfing \cite{ioc}. This series of events followed a documented power struggle between the IOC and the IFSC, as the IFSC sought legitimization for sport climbing through Olympic inclusion \cite{BR19}. Seeking greater levels of funding, professionalization, and notoriety---and having been less than successful with its solo ventures to do so---the IFSC ``allowed the IOC to obtain a degree of power over the sport and traded its autonomy to some extent'' to partially fulfill these aims \cite[p.\ 1684]{BR19}. Although Olympic sport climbing remains governed by the IFSC, the IOC has imposed certain organizational standards onto sport climbing, which have resulted in tensions of ethics and questions of the IFSC's organizational power in the broader climbing community \cite{BR19}.\footnote{For a comprehensive history of international competitive sport climbing and the organizational tension between the IFSC and the IOC, as well as the implications of the Olympic Agenda 2020, please see both Bautev and Robinson \cite{BR19} and  Thorpe and Wheaton \cite{TW19}.} One of these frequently-questioned standards is the development of the IFSC sport climbing combined format, which was influenced by an IOC recommendation (Degun, as cited in \cite{BR19}).

Currently, IFSC sport climbing consists of three individual disciplines that are then scored in a combined format. These disciplines include speed climbing, bouldering, and lead climbing. Speed climbing is a race-format event in which two athletes compete for the fastest time on a 15m fixed route \cite{rules2021}. Speed climbing rewards the fastest time but  also rewards precision, as false starts or falls are automatic losses in finals \cite{rules2021}. Bouldering consists of athletes ``solving'' either four (in qualifications) or three (in finals) boulder ``problems'' roughly 4.5m in height, where they are rewarded for completing these climbs in the least number of attempts, with a halfway ``zone hold'' to further separate the field through partial attempts \cite{rules2021}. Climbing here is done unroped, one at a time, and without safety equipment save for padded mats. Bouldering rewards athleticism, strength, quick-thinking, and adaptability. Lead climbing is perhaps the event that most outwardly resembles traditional rock climbing. Athletes compete one at a time, roped, on a 15m wall on a progressively-difficult course \cite{rules2021}. The athlete to climb the highest wins the event; if two climbers reach the same point on the wall, the quickest athlete is rewarded (with a total available climbing time of six minutes) \cite{rules2021}. Lead climbing rewards endurance and precision---once an athlete falls, their attempt is over. In each discipline, athletes are assigned a ranked score \cite{rules2021}.

Despite the wishes of the IFSC, the IOC only granted one metal per gender to sport climbing, citing limitations due to crowding if each discipline were to have its own medal \cite{BD}. This format brings each athlete's performance in all three disciplines together under one score---the combined format. The combined format remains controversial among the climbing community, and it was publicly denounced by numerous high-profile climbers, including eventual sport climbing Olympians Adam Ondra and Jakob Schubert \cite{BR19,BD}. The controversy surrounds both the clustering of the three disciplines into one event, as well as the specific inclusion of speed climbing in the event in general.  In particular, speed climbing is noted as the ``outlier discipline'' and the ``proverbial wrench in the whole system,'' as ``it's a discipline of climbing that resembles very little traits of outdoor climbing'' \cite{BD}.

And yet the inclusion of speed climbing was necessary in order to offer all sport climbing athletes equal possibility of participation in the Olympic Games while also working with the single-medal quota imposed by the IOC \cite{BD}. Whether or not speed climbing is a ``legitimate'' form of climbing is not relevant to the conversation of scoring, save from the cascade effect of scoring outcomes that came from its inclusion in the combined format.
The interesting part of the question of whether or not speed climbing should be part of the combined format---or what happens to the scoring when it is---lies with the fact that traditionally, speed climbing has been dominated by athletes who are notably uncompetitive in the other two disciplines \cite{BD}. While there are similar discrepancies among individual athletes' skills in bouldering and lead climbing, there is more athlete crossover in those two disciplines than into the speed discipline. However, the announcement of the combined format resulted in many sport climbers taking seriously all three disciplines, and some very respectable all-around climbers have since emerged in both the women's and men's fields.

The inclusion of speed climbing in the combined format likely informed the choice of the current multiplicative scoring system, though finding official documentation of this process proves to be difficult \cite{plastic}. The current scoring for the combined format uses a multiplicative system that takes into account a climber's ranked score (through ``overall ranking points'') in each discipline \cite{rules2018}. This multiplicative system was introduced in April of 2018 following an IFSC Rules Modification in advance of the debut of the Combined format at the 2018 IFSC Climbing World Championships in Innsbruck, Austria \cite{BD,rules2018}). This scoring system has been noted as confusing and anti-climactic, and it was even publicly misunderstood as favouring all-around athletes at first glance \cite{BD,epic}. Despite the possible perception of a combined event being structured to reward consistency across disciplines, the multiplicative system works such that climbers are actually rewarded for being dominant in one discipline, as opposed to being all-around athletes. As Black Diamond \cite{BD} explains, ``[e]ven just one first place finish significantly increases your chances of having a low score, which [favours] the best climbers.''

Our paper herein deals with the possibilities of alternative scorings for the combined format. These alternative scorings include the currently-used multiplicative scoring method, the additive scoring method used prior to 2018, and a new approach that we term the \emph{square root method}.

\subsection{Our Contributions}

The remainder of this paper is organized as follows. In Section \ref{multi.sec}, we provide a brief summary of multi-event scoring methods, including the additive ranking-based scoring systems that are the main subject of this paper. Section \ref{2020.sec} analyzes the results of the 2020 Olympics, comparing the product-based ranking that was used there to the more traditional sum-based rankings. In Section \ref{improved.sec}, we introduce and analyze a square-root based ranking system, which can be viewed as a compromise between the two other scoring methods. We also discuss how additive ranking-based scoring methods can conveniently be implemented using precomputed \emph{scoring tables}. Finally, Section \ref{summary.sec} summarizes our findings and conclusions.

\section{Multi-event Scoring Methods}
\label{multi.sec}

There are many sporting events where the final standings are based on multiple disciplines or on multiple stages of the same discipline. For example, the men's decathlon consists of ten different track and field events. A diving competition may consist of five or six dives (each dive is termed a \emph{round}). It is very common to derive a score for each round of the competition and then compute the sum of each competitor's scores to obtain the final standings. The scores of each round are numerical values that typically fall within some prespecified range.  

There are a much smaller number of sports where the outcomes of each discipline or round are used only to determine a \emph{ranking} of the competitors, and the final outcome only depends on these rankings (sometimes the ranking are called \emph{ordinals}). In this paper, we will refer to such a scoring system as an \emph{ranking-based scoring system}. 

Sport climbing was introduced as an Olympic sport at the 2020 Games (which were held in 2021 due to the Covid-19 pandemic).
As discussed in  Section \ref{intro.sec}, sport climbing consists of three disciplines: speed, bouldering and lead. Each climber competes in all three disciplines, and the final rankings are  determined by \emph{multiplying} the placements in each discipline 
(the lowest score determines the ultimate winner).
In other sports using ranking-based scoring systems, it is more common to compute the \emph{sum} of the rankings in the component disciplines.

We now present a general mathematical description of certain ranking-based scoring systems based on an additive function.
Suppose a sporting event consists of $s$ stages and there are $n$ competitors, each of whom competes in each stage.  For  $1 \leq j \leq n$ and $1 \leq i \leq s$, let $r_{j,i}$ denote the \emph{rank} of the $j$th player in the $i$th stage (a rank is an integer between $1$ and $n$). The \emph{rank vector} for player $j$ is the $s$-tuple $\mathbf{r}_j = (r_{j,1} , \dots , r_{j,s})$. For convenience, and to simplify the discussion, we assume that there are no ties in any stage, so each $n$-tuple $(r_{1,i} , \dots , r_{n,i})$ is a permutation of $\{1, \dots , n\}$, for $1 \leq i \leq s$. 

Let $f: \{1, \dots , n\} \rightarrow \mathbb{R}^+ \cup \{0\}$ be a monotone increasing function; we call $f$  the \emph{score function}.\footnote{The score function is  monotone increasing because we want an $i$th-place finish in any give stage to score less than than an $(i+1)$st-place finish in the same stage (a lower score is better).}
The most common choice for a score function  is the linear function $f(j)= j$ for $1 \leq j \leq n$. The \emph{$f$-score} of player $j$ is the quantity
\[ \mathsf{score}_j = \sum_{i=1}^s f(r_{j,i}). \]
The final ranking of the $n$ competitors is determined by sorting the list of values $\mathsf{score}_j$ in \emph{increasing} order. We note that there may also need to be a tie-breaking mechanism, if
 $\mathsf{score}_j = \mathsf{score}_k$ for some $j \neq k$.

The above definition gives equal weight to each stage. A generalization is to specify a 
\emph{weight vector} $(w_1, \dots , w_s)$ and define the final scores to be
\[ \mathsf{score}_j = \sum_{i=1}^s w_i f(r_{j,i}). \]
We will call this a \emph{weighted score}.
Observe that we obtain the original formula if $w_1 = \dots = w_s = 1$; we could call such a score an \emph{unweighted score}. 

\begin{Example}
{\rm Prior to 2004, figure skating used a weighted additive ranking-based scoring system. Each figure skating competition consisted of a \emph{short program} and a \emph{long program}. The score function was the linear function  $f(j)= j$ (for $j = 1,24$), but the long program received twice the weight of the short program. The rank in the long program was used to break any ties that arose.

One consequence of this scoring system is that any of the top three skaters in the short program
could win the competition by subsequently winning the long program. For example, the total score of a skater who finished third in the short program and first in the long program is $1 \times 3 + 2 \times 1 = 5$. The best total score any other skater could obtain would be $1 \times 1 + 2 \times 2 = 5$; however, in this  case, the skater who won the long program would then be declared the winner.
}
\end{Example}

\begin{Example}
\label{sailing.exam}
{\rm A sailing regatta typically consists of a series of races using an additive ranking-based scoring system. The score function is often, but not always, the linear function $f(j)= j$. In the 1968 Olympics, the scoring function was defined as follows: $f(1) = 0$, $f(2)= 3$, $f(3) = 5.7$, $f(4) = 8$, $f(5) = 10$, $f(6) = 11.7$, and $f(j) = j+6$ if $j \geq 7$. From this scoring system, it can be inferred that a first- and third-place finish in two races is considered to be better than two second-place finishes, because $0 + 5.7 < 2 \times 3$.}
\end{Example}

\begin{Example}
{\rm
The William Lowell Putnam Mathematical Competition \cite{Putnam} is an annual written mathematics competition for undergraduate mathematics students in Canada and the U.S. Each student receives a score between $0$ and $120$. This determines a ranking of all the students who took part in the competition. Before 2019, each university could also designate a 3-person team before the competition took place. The team score was obtained by computing the sum of the rankings of the three students in the team.\footnote{Starting in 2019, a new team scoring system  
was used, in which the sums of the \emph{scores} of the team members was used.}}
\end{Example}

We  already mentioned in Section \ref{intro.sec} that sport climbing in the 2020 Olympics used a product-based scoring system. 
The score function is the usual linear function  $f(j)= j$, but a 
player's score is the product of their three scores (or rankings):
\[ \mathsf{score}_j = \prod_{i=1}^3 r_{j,i}. \] However, we can easily see that there is an equivalent ranking function for sport climbing that is just an additive ranking-based scoring system with a \emph{nonlinear} scoring function. 

We have 
\begin{eqnarray*}
\mathsf{score}_j \leq \mathsf{score}_k &\Leftrightarrow & \prod_{i=1}^3 r_{j,i} \leq \prod_{i=1}^3 r_{k,i}\\
& \Leftrightarrow & \ln \left( \prod_{i=1}^3 r_{j,i}\right) \leq \ln \left( \prod_{i=1}^3 r_{k,i} \right)\\
& \Leftrightarrow & \sum_{i=1}^3 \ln  r_{j,i} \leq \sum_{i=1}^3 \ln  r_{k,i}.
\end{eqnarray*}
Thus, if we use a \emph{logarithmic} scoring function, $f(j) = \ln j$, then the resulting 
additive ranking-based scoring system yields the same final rankings as the previously described multiplicative ranking-based scoring system.

For computations, it is probably simplest to compute the product of three rankings as opposed to computing the sum of their logarithms. Furthermore, it most sports announcers on television would probably not be comfortable discussing logarithms. However, we can gain some insight into the properties of the sum vs the product scoring system by recognizing that the product scoring system is just an additive system with a different score function. We will discuss this further in  Section \ref{improved.sec}.

\section{Analysis of Results at the 2020 Olympics}
\label{2020.sec}

Tables \ref{tab3}--\ref{tab6} show two possible sets of outcomes of the sport climbing  preliminaries and finals (men's and women's) at the 2020 Olympics. Note that at the 2020 Olympics, preliminaries were used to reduce the number of competitors from 20 to 8. 
The finals then involved the eight best climbers from the preliminaries.\footnote{It is important to note that in the men's final, the seventh place qualifier, B.\ Mawem, did not compete, due to a torn bicep injury sustained during his last climb of the qualification round. B.\ Mawem was marked Did Not Start (DNS) for the finals round, but according to IOC rules he still finished 8th overall.}

First, we give the official rankings as determined by multiplying the discipline rankings. The second (hypothetical) set of rankings uses the more common method of computing the sum of the discipline rankings. Each triple of discipline rankings consists of the rankings for speed, bouldering and lead (in that order).

\begin{table}
\caption{Sport Climbing Men's Preliminaries Sum vs Product Rankings}
\label{tab3}
\begin{center}
\begin{tabular}{|l|c||r|r||r|r|}
\hline
 Name & Discipline Rankings & Product & Ranking & Sum & Ranking \\ \hline \hline
	M.\ Mawem 	 &	$(3	 ,	1	 ,	11)$	 &	33	& 1 & 15 & 2 \\\hline
		Narasaki 	 &	$(2	 ,	2	 ,	14)$	 &	56	& 2 & 18 & 3\\\hline
		Duffy 	 &	$(6	 ,	5	 ,	2)$	 &	60	& 3 & 13 & 1\\\hline
	Schubert 	 &	$(12	 ,	7	 ,	1)$	 &	84	 & 4 & 20 & 4\\\hline
		Ondra 	 &	$(18	 ,	3	 ,	4)$	 &	216	 & 5& 25 & 6\\\hline
	Gin\'{e}s L\'{o}pez 	 &	$(7	 ,	14	 ,	3)$	 &	294	 & 6 & 24 & 5\\\hline
	B.\ Mawem 	 &	$(1	 ,	18	 ,	20)$	 &	360	& 7  &\textcolor{red}{39} & \textcolor{red}{17}\\\hline
	Coleman 	 &	$(10	 ,	11	 ,	5)$	 &	550	 &8 & 26 & 7\\\hline \hline
	Megos 	 &	$(19	 ,	6	 ,	6)$	 &	684	& 9 & 31 & 9 \text{(tie)} \\\hline
		Chon 	 &	$(5	 ,	10	 ,	16)$	 &	800	 & 10& 31 & 9 \text{(tie)}\\\hline
		Khaibullin 	 &	$(4	 ,	17	 ,	13)$	 &	884	 & 11& 34 & 12\\\hline
	Hojer 	 &	$(11	 ,	9	 ,	9)$	 &	891	 & 12& \textcolor{red}{29} & \textcolor{red}{8}\\\hline
		Rubtsov 	 &	$(16	 ,	4	 ,	15)$	 &	960	 &13& 35 & 13 \text{(tie)}\\\hline
		Pan 	 &	$(20	 ,	8	 ,	7)$	 &	1120	&14 & 35 & 13 \text{(tie)}\\\hline
		Piccolruaz 	 &		$(8	 ,	13	 ,	12)$	 &	1248	& 15 & 33 & 11\\\hline
		Cosser 	 &		$(9	 ,	16	 ,	10)$	 &	1440	& 16 & 35 & 13 \text{(tie)}\\\hline
		McColl 	 &	$(14	 ,	15	,	8)$	 &	1680	& 17 & 37 & 16\\\hline
		Harada 	 &	$(15	 ,	12	 ,	17)$	 &	3060	& 18 & 44 & 18\\\hline
		Fossali 	 &		$(13	 ,	19.5	 ,	18)$	 &	4563	& 19 & 50.5 & 19\\\hline
		O'Halloran 	 &	$(17	 ,	19.5	 ,	19)$	 &	6298.5 & 20 & 55.5 & 20\\\hline
\end{tabular}
\end{center}
\end{table}

\begin{table}
\caption{Sport Climbing Men's Finals Sum vs Product Rankings}
\label{tab4}
\begin{center}
\begin{tabular}{|l|c||r|r||r|r|}
\hline
 Name & Discipline Rankings & Product & Ranking & Sum & Ranking \\ \hline \hline
	Gin\'{e}s L\'{o}pez 	 &	$(1,7,4)$	 &	28	 & 1 & 12 & 2 \text{(tie})\\\hline
		Coleman 	 &	$(6,1 ,	5)$	 &	30	 &2 & 12 & 2 \text{(tie)}\\\hline 
	Schubert 	 &	$(7,5	 ,	1)$	 &	35	 & 3 & 13 & 7\\\hline
		Narasaki 	 &	$(2	 ,	3,6)$	 &	36	& 4 & 11 & 1\\\hline
M.\ Mawem 	 &	$(3	 ,	2,7)$	 &	42	& 5 & 12 & 2 \text{(tie)} \\\hline
		Ondra 	 &	$(4,6,2)$	 &	48	 & 6& 12 & 2 \text{(tie)}\\\hline
		Duffy 	 &	$(5,4,3)$	 &	60	& 7 & 12 & 2 \text{(tie)}\\\hline
\end{tabular}
\end{center}
\end{table}

\begin{table}
\caption{Sport Climbing Women's Preliminaries Sum vs Product Rankings}
\label{tab5}
\begin{center}
\begin{tabular}{|l|c||r|r||r|r|}
\hline
 Name & Discipline Rankings & Product & Ranking & Sum & Ranking \\ \hline \hline
	Garnbret 	 &	$(14,1,4)$	 &	56	& 1 & 19 & 3 \\\hline
		Seo 	 &	$(17,5,1)$	 &	85	& 2 & 23 & 6\\\hline
		Nonaka 	 &	$(4,8,3)$	 &	96	& 3 & 15 & 1\\\hline
	Noguchi 	 &	$(9,3,6)$	 &	162	 & 4 & 18 & 2\\\hline
		Raboutou 	 &	$(12,2,8)$	 &	192	 & 5& 22 & 4 \text{(tie)} \\\hline
	Pilz 	 &	$(11,9,2)$	 &	198	 & 6 & 22 & 4 \text{(tie)}\\\hline
	Miroslaw 	 &	$(1	 ,	20,19)$	 &	380	& 7  & \textcolor{red}{40} & \textcolor{red}{16 \text{(tie)}}\\\hline
	Jaubert 	 &	$(2,13,15)$	 &	390	 &8 & \textcolor{red}{30} & \textcolor{red}{9}\\\hline \hline
	Meshkova 	 &	$(15,6,5)$	 &	450	& 9 & \textcolor{red}{26} & \textcolor{red}{7}  \\\hline
		Coxsey 	 &	$(16,4,13)$	 &	832	 & 10& 33 & 11 \\\hline
		Condie 	 &	$(7,11,11)$	 &	847	 & 11& \textcolor{red}{29} & \textcolor{red}{8}\\\hline
	Song 	 &	$(3,19,18)$	 &	1026	 & 12& 40 & 16 \text{(tie)}\\\hline
		Chanourdie 	 &	$(8,15,9)$	 &	1080 &13& 32 & 10 \\\hline
		Yip 	 &	$(6,16,12)$	 &	1152	&14 & 34 & 12 \text{(tie)} \\\hline
		Rogora 	 &		$(19,7,10)$	 &	1330	& 15 & 36 & 14\\\hline
		Klingler 	 &		$(10,10,14)$	 &	1400	& 16 & 34 & 12 \text{(tie)}\\\hline
		Kaplina 	 &	$(5,18,17)$	 &	1530	& 17 & 40 & 16 \text{(tie)}\\\hline
		Krampl 	 &	$(18,14,7)$	 &	1764	& 18 & 39 & 15\\\hline
		MacKenzie 	 &		$(13	 ,	12,16)$	 &	2496	& 19 & 41 & 19\\\hline
		Sterkenburg 	 &	$(20,17,20)$	 &	6800 & 20 & 57  & 20\\\hline
\end{tabular}
\end{center}
\end{table}

\begin{table}
\caption{Sport Climbing Women's Finals Sum vs Product Rankings}
\label{tab6}
\begin{center}
\begin{tabular}{|l|c||r|r||r|r|}
\hline
 Name & Discipline Rankings & Product & Ranking & Sum & Ranking \\ \hline \hline
	Garnbret	 &	$(5,1,1)$	 &	5	 & 1 & 7 & 1 \\\hline
		Nonaka 	 &	$(3,3,5)$	 &	45	 &2 & 11 & 2 \\\hline 
	Noguchi 	 &	$(4,4,4)$	 &	64	 & 3 & 12 & 3\\\hline
		Miroslaw 	 &	$(1,8,8)$	 &	64	& 4 & 17 &7 \text{(tie)}\\\hline
Raboutou 	 &	$(7,2,6)$	 &	84	& 5     & 15 & 5 \text{(tie)} \\\hline
		Jaubert 	 &	$(2,6,7)$	 &	84	 & 6& 15 & 5 \text{(tie)}\\\hline
		Pilz 	 &	$(6,5,3)$	 &	90	& 7     & 14 & 4\\\hline
		Seo 	 &	$(8,7,2)$	 &	112	& 8     & 17 & 7 \text{(tie)}\\\hline
\end{tabular}
\end{center}
\end{table}

It could be argued that the choice to multiply discipline rankings was made by the IFSC because under the former additive system, no speed specialists would qualify for the finals \cite{plastic}. This is because many of the top speed specialists are not as competitive in the more technical disciplines of bouldering and lead climbing, something that is most evident with men's competitor B. Mawem, and women's competitor Miroslaw.  Thus, the modified (multiplicative) system was employed, in the hope that this would lead to some speed specialists qualifying for the final round.

The main effect of multiplying rankings is that it places a very large premium on finishing first in a discipline. (For example, if a first place finish is replaced by a second place finish, then the overall score is doubled.) A competitor who finishes first in a discipline is very likely to qualify for the finals, even if their finish is close to the bottom in the other two disciplines. So, in this respect, the modified scoring system achieved its desired goal. Unfortunately, at the same time, it could be argued that multiplying rankings tends to undervalue to a certain extent an all-around climber who is quite good but not outstanding in all three disciplines. This seems directly contrary to what should be the purpose of a combined event.

Tables \ref{tab3}--\ref{tab6} illustrate how the outcomes would have differed in the 2020 Olympics  in the two scoring systems. 
It should be emphasized that many of the final rankings are similar or roughly similar in 
both scoring systems. But examining the differences and identifying the outliers is interesting and instructive, particularly if we wish to develop a scoring system more reflective of the aims of a combined format (i.e., finding the best overall athlete).

First, we look at the men's preliminary round. Recall that the purpose of the preliminary round is to reduce the size of the field from 20 to 8. From Table \ref{tab3}, we see that the main difference between the results of two scoring methods is that B.\ Mawem would have been replaced by Hojer if a sum-based scoring system had been used.\footnote{B.\ Mawem did not compete in the finals due to injury, so there were only seven climbers in the final.} B.\ Mawem won the speed discipline and finished 18th and 20th in the other two disciplines. Thus, he ended up in the top eight according to the product score, but he would have finished 17th out of 20 if the sum score had been used instead. Had the sum score been used, B.\ Mawem would have been replaced by Hojer, who had three ``middle-of-the-pack'' finishes, namely, 11th, 9th and 9th. We think that despite the decisions of the IFSC, many people would find it problematic that someone who combined a first place finish with two very low finishes should advance to the final in a combined event, while someone who is competent but not outstanding an all three disciplines is passed over. 

When we turn to the men's finals, we find that the three top placements were obtained by the three competitors who won one of the disciplines. Again, the premium for finishing first in a discipline outweighs significantly poorer placements in the other disciplines, something that we see even with the Gold medal finisher Gin\'{e}s L\'{o}pez, who scored first in speed climbing yet last in the bouldering round. Gin\'{e}s L\'{o}pez's first-place finish in speed was a particular boon for him as he is known as a lead climber, the discipline in which he finished fourth \cite{plastic}. 

If we instead computed the sum of the three rankings, we see that the fourth-place finisher (Narasaki) would have won. Narasaki is notable as a cross-disciplinary athlete, and perhaps the most accomplished non-speed specialist in speed climbing; he developed and popularized a unique way of moving through the speed route, deemed the ``Tomoa skip'' \cite{Samet}. (Ironically, it was this move that he fumbled in his race against Gin\'{e}s L\'{o}pez, ultimately resulting in his second-place speed finish.) Narasaki would have been followed by five climbers who tied for second place (of course a tie-breaking mechanism would be employed to separate the finishes of these five climbers, e.g., a count-back to their qualification standings). The third-place finisher (Schubert) would have finished last if the final ranking had been based on the sum of the rankings.

In the women's competition, similar differences can be found between the two scoring systems. In the preliminary round, the 7th and 8th finishers (Miroslaw and Jaubert) both combined one high finish
(first or second) with two below average finishes, but this enabled them to qualify for the finals. The seventh place finisher, Miroslaw, won the speed event but finished 19th and 20th in the other two disciplines. She would have finished in a three-way tie for 16th if the sum scoring system had been used. The two competitors (namely, Meshkova and Condie) who would have replaced Miroslaw and Jaubert (had the sum scoring system been used) both had more ``uniform'' finishes in the three disciplines.

The three medalists in women's final would have been the same under both scoring systems. Garnbret was the  favourite to win the gold medal, and indeed she was exceptionally dominant with her discipline rankings of $(5,1,1)$. Nonaka and Noguchi were both particularly consistent cross-discipline, which was also not unexpected. Indeed, just before the Games, Nonaka became one of the first women's non-speed specialists to podium at an IFSC World Cup speed event \cite{Walker}. The most significant discrepancy is that Miroslaw would have finished in a tie for last place under the sum system, instead of finishing fourth (she combined a first-place finish with two last-place finishes in the finals). As a true speed specialist, Miroslaw is perhaps the most overt example of the speed vs bouldering/lead inconsistency. In her final speed run, Miroslaw set a new women's speed world record with a time of 6.84 seconds. However, in the bouldering final round shortly after, Miroslaw was unable to score a single zone (finishing with a score of 0), and she fell off the lead route at hold 9+, only a quarter of the progress of first-place finisher Garnbret. In another notable ranking shift, Pils would have moved up from 7th place to 4th place.

The women's finals also included two two-way ties, one for third and fourth place, and one for fifth and sixth place. A two-way tie was broken by comparing the head-to-head finishes; the climber who won two out of three of these was ranked higher \cite{rules2021}. Thus Raboutou was ranked above Jaubert and Noguchi was ranked above Miroslaw. As a result, Noguchi won the bronze medal. It is interesting to compare 
the rankings of Noguchi and Miroslaw: Noguchi's rankings were $(4,4,4)$ while Miroslaw's were $(1,8,8)$. In this particular case, the tie-breaking mechanism favoured the climber with three equal finishes over the climber with one first-place and two last place finishes. Arguably this is a reasonable result, but it is contrary to the apparent goal of the product system to give preference to first-place results.

\section{An Alternative Ranking-Based Scoring Method}
\label{improved.sec}

As can be seen from the analysis of the data sets that we carried out in Section \ref{2020.sec}, the product scoring system enabled some speed specialists to achieve much higher finishes (indeed, any climber who finishes first in one discipline and low in the other two disciplines would benefit greatly). But we question whether this is completely fair in the context of a combined event. On the other hand, the sum scoring system tends to undervalue first place finishes. For example a first and third place finish is treated the same as two second-place finishes (see Example \ref{sailing.exam} for a different scoring method in the setting of sailing competitions that intentionally avoids this scenario). Thus we think it would be useful to consider an alternate scoring system.

To further illustrate, let us consider when a first- and last-place is ranking in two disciplines is equivalent to two ``similar'' rankings. If we use the product scoring system, we see that a first- and 20th-place finish is equivalent to a fourth- and fifth-place finish, because $1 \times 20 = 4 \times 5$. On the other hand, in the sum scoring system, a first- and 20th-place finish is equivalent to a 10th- and 11th-place finish, because $1 + 20 = 10 + 11$. It might be preferable to have a scoring system that achieves more of a compromise, e.g., one in which a first- and 20th-place finish is (roughly) equivalent to a 6th- and 7th-place finish, or a 7th- and 8th-place finish. 

For the time being, it will be useful to consider additive systems, so we will speak in terms of the logarithmic scoring function instead of the product system (recall that they lead to identical rankings).  Given the drawbacks of the linear and logarithmic scoring functions, we could instead consider a scoring function that is between them. Basically, we would seek a concave function, but one that is less concave than the logarithm function.\footnote{A real-valued function is  \emph{concave} if the line segment between any two points on the graph of the function lies below  the graph between the two points.}

\begin{figure}[t]
\begin{center}


\includegraphics[width=3.0in]{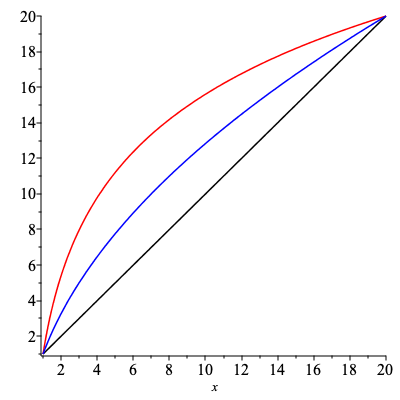}
\end{center}

\[
\begin{array}{ll}
\text{black} & g_1(j) \text{ (linear) }\\
\text{blue} & g_2(j) \text{ (square root) }\\
\text{red} & g_3(j) \text{ (logarithmic) }\\
\end{array}
\]

\caption{Three possible scoring functions}
\label{scoringfunctions.fig}
\end{figure}

The function $f(j) =\sqrt {j}$ is a reasonable choice. (More generally, we could employ a function of the form $f(j) = j^f$, where $0<f < 1$ is a fixed real number.)  We  compute
\begin{eqnarray*}
\sqrt{1} + \sqrt{20} &=& 5.472\\
\sqrt{7} + \sqrt{8} &=& 5.474.
\end{eqnarray*}
Thus, this square root scoring function treats a 7th- and 8th-place finish as basically equivalent to a first- and 20th-place finish, as we suggested above. 

In Figure \ref{scoringfunctions.fig}, we illustrate how a square root scoring function lies between a linear and a logarithmic scoring function. We want to compare the three functions $f_1(j) = j$, 
$f_2(j) = \sqrt{j}$ and $f_3(j) = \ln j$. To obtain a nice visual comparison, we adjust the three functions via  affine transformations so that 
$f_1(1) = f_2(1) = f_3(1) = 1$ and $f_1(20) = f_2(20) = f_3(20) = 20$. The affine transformations do not affect any resulting rankings. So we are actually comparing the following three functions in Figure \ref{scoringfunctions.fig}:
\begin{eqnarray*}
g_1(j) &=& j\\
g_2(j) &=& \frac{\sqrt{20}- 20}{\sqrt{20}- 1} + \frac{19}{\sqrt{20}- 1} \sqrt{j}\\
g_3(j) &=& 1 + \frac{19}{\ln 20}\ln j.
\end{eqnarray*}

It is interesting to compare the rankings obtained from the square root scoring function to the rankings we obtained previously in Section \ref{2020.sec}. For the men's preliminaries, the square root ranking would have qualified Megos while demoting B.\ Mawem.
Hojer would have moved up, but only to 9th place instead of 8th place. (See Table \ref{tab7} for the complete rankings.)  For the men's finals, the square root ranking is essentially the same as the sum ranking: Narasaki would be ranked first, followed by Gin\'{e}s L\'{o}pez and Coleman (see Table \ref{tab8}).

The women's results are found in Tables \ref{tab9} and \ref{tab10}. In the preliminaries, Meshkova would have qualified instead of  Miroslaw, as with the sum ranking. However, the second swap in the sum ranking (Jaubert for Condie) does not occur in the square root ranking. Thus, the square root ranking is a compromise between the sum and product rankings. For the finals, the three medal winners are identical in all three scoring methods considered.

\begin{table}
\caption{Sport Climbing Men's Preliminaries Rankings Including Square Root Scores}
\label{tab7} 
\begin{center}
\begin{tabular}{|l|c|c||r|r|r|}
\hline
 Name & Discipline Rankings &  \multicolumn{1}{|c||}{$\sqrt{\quad}$ score}  & \multicolumn{3}{|c|}{Overall Rankings}\\
&  &   &  Product  & Sum   & $\sqrt{\quad}$ \\ \hline \hline
 	M.\ Mawem 	 &	$(3	 ,	1	 ,	11)$ &	6.049& 1 &   2 & 1\\\hline
		Narasaki 	 &	$(2	 ,	2	 ,	14)$	& 	6.570	& 2 &   3& 3\\\hline
		Duffy 	 &	$(6	 ,	5	 ,	2)$	 &		6.100 & 3 &   1& 2\\\hline
	Schubert 	 &	$(12	 ,	7	 ,	1)$	 &	7.110 &	  4 &   4& 4\\\hline
		Ondra 	 &	$(18	 ,	3	 ,	4)$	 &7.975	&	  5&   6& 5\\\hline
	Gin\'{e}s L\'{o}pez	 &	$(7	 ,	14	 ,	3)$	 &8.119	&	  6 &   5& 6\\\hline
	B.\ Mawem 	 &	$(1	 ,	18	 ,	20)$	 &	9.715&	 7   & \textcolor{red}{17}& \textcolor{red}{11}\\\hline
	Coleman 	 &	$(10	 ,	11	 ,	5)$	 &	8.715&	 8 &   7&   7\\\hline \hline
	Megos 	 &	$(19	 ,	6	 ,	6)$	 &	9.258	& 9   & 9 \text{(tie)} &   \textcolor{red}{8}\\\hline
		Chon 	 &	$(5	 ,	10	 ,	16)$	 &	9.398	&  10&   9 \text{(tie)}& 10\\\hline
		Khaibullin 	 &	$(4	 ,	17	 ,	13)$	 & 9.729	&	  11&   12& 12\\\hline
	Hojer 	 &	$(11	 ,	9	 ,	9)$	 &	9.317	&  12&   \textcolor{red}{8}& 9\\\hline
		Rubtsov 	 &	$(16	 ,	4	 ,	15)$	 & 9.873	&	 13&   13 \text{(tie)}& 13\\\hline
		Pan 	 &	$(20	 ,	8	 ,	7)$	 &	9.946	&14   & 13 \text{(tie)}& 15\\\hline
		Piccolruaz 	 &		$(8	 ,	13	 ,	12)$&	9.898 &		 15 &   11& 14\\\hline
		Cosser 	 &		$(9	 ,	16	 ,	10)$	 &	10.162 &	 16   & 13 \text{(tie)}& 16\\\hline
		McColl 	 &	$(14	 ,	15	,	8)$	 &	10.443	& 17 &   16& 17\\\hline
		Harada 	 &	$(15	 ,	12	 ,	17)$	 &	11.460	& 18 &   18& 18\\\hline
		Fossali 	 &		$(13	 ,	19.5	 ,	18)$	& 	12.264	& 19 &   19& 19\\\hline
		O'Halloran 	 &	$(17	 ,	19.5	 ,	19)$	 &	12.898 &  20 &   20& 20\\\hline
\end{tabular}
\end{center}
\end{table}

\begin{table}
\caption{Sport Climbing Men's Finals Rankings Including Square Root Scores}
\label{tab8}
\begin{center}
\begin{tabular}{|l|c|c||r|r|r|}
\hline
 Name & Discipline Rankings &  \multicolumn{1}{|c||}{$\sqrt{\quad}$ score}  & \multicolumn{3}{|c|}{Overall Rankings}\\
&  &   &  Product  & Sum   & $\sqrt{\quad}$ \\ \hline \hline
 Gin\'{e}s L\'{o}pez		 &	$(1,7,4)$ &	5.646& 1 &   2 \text{(tie)} & 2\\\hline
		Coleman 	 &	$(6,1,5)$	& 	5.686	& 2 &   2 \text{(tie)}& 3\\\hline
	Schubert	 &	$(7,5,1)$	 &		5.882 & 3 &   7& 6\\\hline
	 	Narasaki &	$(2,3,6)$	 &	5.596 &	  4 &   1& 1\\\hline
M.\ Mawem 			 &	$(3,2,7)$	 &5.792	&	  5&   2 \text{(tie)}& 4\\\hline
	Ondra 	 &	$(4,6,2)$	 &5.864	&	  6 &   2 \text{(tie)}& 5\\\hline
	Duffy 		 &	$(5,4,3)$	 &5.968	&	  7 &   2 \text{(tie)} & 7\\\hline
\end{tabular}
\end{center}
\end{table}

\begin{table}
\caption{Sport Climbing Women's Preliminaries Rankings Including Square Root Scores}
\label{tab9}
\begin{center}
\begin{tabular}{|l|c|c||r|r|r|}
\hline
 Name & Discipline Rankings &  \multicolumn{1}{|c||}{$\sqrt{\quad}$ score}  & \multicolumn{3}{|c|}{Overall Rankings}\\
&  &   &  Product  & Sum   & $\sqrt{\quad}$ \\ \hline \hline
 	Garnbret 	 &	$(14,1,4)$ &	6.742 & 1 &   3 & 2\\\hline
		Seo 	 &	$(17,5,1)$		& 	7.359	& 2 &   6& 4\\\hline
		Nonaka 	 &	$(4,8,3)$		 &		6.560& 3 &   1& 1\\\hline
	Noguchi 	 &	$(9,3,6)$	 &	7.182 &	  4 &   2& 3\\\hline
		Raboutou 	 &	$(12,2,8)$	 &7.707	&	  5&   4 \text{(tie)}& 5\\\hline
	Pilz 	 &	$(11,9,2)$	 &7.731	&	  6 &   4 \text{(tie)}& 6\\\hline
	Miroslaw 	 &	$(1	 ,	20,19)$	 &	9.831&	 7   & \textcolor{red}{16}& \textcolor{red}{12}\\\hline
	Jaubert 	 &	$(2,13,15)$		 &	8.893&	 8 &   \textcolor{red}{9} &   8\\\hline \hline
	Meshkova 	 &	$(15,6,5)$	 &	8.559	& 9   & \textcolor{red}{7} &   \textcolor{red}{7}\\\hline
		Coxsey 	 &	$(16,4,13)$		 &	9.606	&  10&   11 & 10\\\hline
		Condie 	 &	$(7,11,11)$		 & 9.279	&	  11&   \textcolor{red}{8} & 9\\\hline
	Song 	 &	$(3,19,18)$		 &	10.334	&  12& 16 \text{(tie)}  & 16\\\hline
		Chanourdie 	 &	$(8,15,9)$		 & 9.701	&	 13 &   10 & 11\\\hline
		Yip 	 &	$(6,16,12)$	 &	9.914	&14   & 12 \text{(tie)} & 13\\\hline
		Rogora 	 &		$(19,7,10)$	 &	10.167 &		 15 &   14 & 15\\\hline
		Klingler 	 &		$(10,10,14)$	 &	10.066 &	 16   & 12 \text{(tie)}& 14\\\hline
		Kaplina 	 &	$(5,18,17)$	 	 &	10.602	& 17 &   16 \text{(tie)}  & 17\\\hline
		Krampl 	 &	$(18,14,7)$		 &	10.630	& 18 &   15 & 18\\\hline
		MacKenzie 	 &		$(13	 ,	12,16)$		& 	11.070	& 19 &   19& 19\\\hline
		Sterkenburg 	 &	$(20,17,20)$		 &	13.067 &  20 &   20& 20\\\hline
\end{tabular}
\end{center}
\end{table}

\begin{table}
\caption{Sport Climbing Women's Finals Rankings Including Square Root Scores}
\label{tab10}
\begin{center}
\begin{tabular}{|l|c|c||r|r|r|}
\hline
 Name & Discipline Rankings &  \multicolumn{1}{|c||}{$\sqrt{\quad}$ score}  & \multicolumn{3}{|c|}{Overall Rankings}\\
&  &   &  Product  & Sum   & $\sqrt{\quad}$ \\ \hline \hline
 Garnbret	 &	$(5,1,1)$	 &	4.236 & 1 &   1 & 1\\\hline
		Nonaka 	 &	$(3,3,5)$ & 	5.700	& 2 &   2& 2\\\hline
	Noguchi 	 &	$(4,4,4)$	  &		6.000 & 3 &   3& 3\\\hline
	 	Miroslaw 	 &	$(1,8,8)$	 &	6.657 &	  4 &   7 \text{(tie)} & 7\\\hline
Raboutou 	 &	$(7,2,6)$ &6.509 (tie)	&	  5&   5 \text{(tie)}& 5 \text{(tie)}\\\hline
	Jaubert 	 &	$(2,6,7)$		 &6.509 (tie)	&	  6 &  5 \text{(tie)}& 5 \text{(tie)}\\\hline
	Pilz 	 &	$(6,5,3)$	 & 6.418	&	  7 &   4 & 4\\\hline
	Seo 	 &	$(8,7,2)$		 &6.888	&	  8 &   7 \text{(tie)} & 8\\\hline
\end{tabular}
\end{center}
\end{table}

\subsection{Complexity of the Scoring Methods}

The three scoring methods we have analyzed are similar in that they can all be viewed as additive ranking-based systems. The only difference is that they employ different scoring functions. Obviously the usual sum-based system is the simplest to understand. As we pointed out, the product-based system is equivalent to computing the sum of the logarithms of the rankings, and we have proposed a new scoring system based on computing the sum of the square roots of the rankings. 

Viewers and commentators who are not mathematically inclined might not be comfortable discussing square roots and logarithms. However, it simple to generate a \emph{scoring table} which lists the points awarded for each ranking in a discipline (e.g., rankings of 1--20 in the preliminaries and 1--8 in the finals). To avoid having to deal with fractions, the relevant logarithms or square roots could be multiplied by 100 or 1000, say, and then rounded to the nearest integer. (This of course would not affect the rankings obtained from these scores.) It should be noted that using a score table is common in other athletic events, e.g., the decathlon, where there are ten different ``performance tables,'' one for each event. The decathlon performance tables convert a time or distance into a numerical score for that event.

Two possible scoring tables are listed in Table \ref{tab11}. We have used the function $100 \ln n$ for the  logarithm-based scores (which yield rankings equivalent to the product-based scoring method) and the function
$100 \sqrt{n} - 100$ for the the square-root based scores. 
These logarithm-based scores  range from $0$ to $347$, while the square-root based scores range from $0$ to $300$, which seem to be a reasonable range of possible values. Of course, these scoring tables could be adjusted to any desired range by applying a suitable affine transformation, which would preserve any rankings obtained using them.

\begin{table}[t]
\caption{Two Possible Scoring Tables}
\label{tab11}
\begin{center}
\begin{tabular}{|r|c|c|}
\hline
 Ranking & Square root-based score & Logarithm-based score  \\ \hline \hline
                           1 &  0 &  0\\ \hline 
                  2 &  41 &  69\\ \hline
                  3 &  73 &  110\\ \hline
                  4 &  100 &  139\\ \hline
                  5 &  124 &  161\\ \hline
                  6 &  145 &  179 \\ \hline
                  7 &  165 &  195\\ \hline
                  8 &  183 &  208\\ \hline
                  9 &  200 &  220\\ \hline
                  10 &  216  &  230\\ \hline
                  11 &  232 &  240\\ \hline
                  12 &  246  &  248\\ \hline
                  13 &  261 &  256 \\ \hline
                  14 &  274  &  264\\ \hline
                  15 &  287  &  271\\ \hline
                  16 &  300  &  277\\ \hline
                  17 &  312  &  283\\ \hline
                  18 &  324  &  289\\ \hline
                  19 &  336 &  294\\ \hline
                  20 &  347 &  300\\ \hline
\end{tabular}
\end{center}
\end{table}

\section{Summary and Conclusion}
\label{summary.sec}

In a follow-up analysis of the 2020 Games by Plastic Weekly, host Tyler Norton expressed the opinion that one of the downfalls of a product-based scoring system was that the mental load of constantly calculating standings eclipsed the performances of many climbers \cite{plastic}. Although the dynamic nature of the multiplicative system resulted in the dramatic shifting of the men's podium based on Schubert's final lead climb, ultimately this took away from the ``raw climbing experience,'' making it ``less about the climbing'' \cite{plastic}. While we empathize with the inherent tensions in and complications of quantifying and rationalizing rock climbing in general, we don't believe that \emph{all} possible sport climbing scoring formats would be equally as distracting as the multiplicative format. With an appropriate scoring system, competitive sport climbing can still be ``about the climbing.''

There are also numerous comments that could be made about the intricacies of each discipline's scoring formats, including how the effects of athlete injury (B. Mawem), speed false starts (Duffy) and slips (Narasaki), and unexpected bouldering performances (Ondra) affected the final rankings, particularly in the men's event.\footnote{For a detailed play-by-play of both the men's and women's finals, please see Climber News \url{https://www.climbernews.com/mens-olympic-climbing-final-results/} and \url{https://www.climbernews.com/womens-olympic-climbing-final-results/}.} Another factor to consider is not athlete performance, but the effect of route-setting (the design and construction of the climbing problems and routes), especially in the bouldering rounds with reference to what has been called ``the Janja problem'' (i.e., that Garnbret so far exceeds the field in bouldering that building bouldering problems that achieve decent separation is difficult. We see this in the women's finals, where Garnbret topped two of the three problems and no one else topped a single problem.) We have intentionally limited the scope of our paper, and thus we do not discuss the effects of the 
individual-discipline scoring and competition rules, though there are likely interesting conversations to pick up about the differences in speed scoring between qualifications (best time) and finals (head-to-head knockout format), as well as bouldering (four boulder problems in qualifiers, three boulder problems in finals). 

The men's final placements have been the subject of much public scrutiny, and indeed much of the conversation surrounding alternative scoring formats post-Olympics was oriented toward trying to manufacture a podium that was more ``publicly acceptable'' than the actual final results \cite{plastic}. (This is not the case with the women's finals, which were widely considered to be an accurate reflection of the field.) It is important to clarify that we are not attempting to add to this conversation to detract from the accomplishments of the winners, but to speak to the disconnect between event aims and goals (i.e., a combined ``overall'' event) and scoring (i.e., scoring that rewards outstanding performance in one discipline). 

Finally, our recommendation for a square root-based scoring method is primarily a theoretical exercise, as the combined event as structured at the Tokyo 2020 Games will not be held again at the Paris 2024 Games \cite{plastic}. Instead, the IOC has granted an additional medal to each gender, and the IFSC has decided to run a speed-only event, and a second event combining bouldering and lead climbing \cite{Paris}. While this does not completely remove the problem of calculating overall scoring across two disciplines, there is much more traditional overlap between bouldering and lead climbing than between speed and either of the other two disciplines. In part, this Paris 2024 two-event format should produce better and more interesting results in both speed climbing and in bouldering/lead climbing. Nevertheless, we wanted to attempt to respond to John Burgman's claim of ``I don't know if anybody has thought of a better system yet,'' with a possible alternative scoring system (namely, the square-root based system) that addresses the problems of  over-valuing single-discipline wins  and under-valuing cross-discipline consistency \cite{plastic}.

\section*{Acknowledgement}

We would like to thank Bill Martin for helpful discussions.


\begin{thebibliography}{XX}


\bibitem{Adie}
R.\ Adie. IFSC Proposes Sport Climbing For Toyko 2020 Olympic Games. 
The British Mountaineering Council, August 8, 2015. \url{https://www.thebmc.co.uk/ifsc-proposes-sport-climbing-for-tokyo-2020}

\bibitem{BR19}
M.\ Bautev and L.\ Robinson. Organizational evolution and the Olympic Games: the case of sport climbing. \emph{Sport in Society} {\bf 22(10)} (2019), 1674--1690.

\bibitem{BD}
Black Diamond. Ondra takes second at the Combined Finals debut during the 2018 World Championship. Black Diamond Equipment. \url{https://www.blackdiamondequipment.com/fr_eu/ stories/experience-story-olympic-format-explained/?cid=olympic-format- explained&utm_campaign=olympicformatexplained&utm_medium=experience&utm_source=in stagram_mktg}


\bibitem{epic}
EpicTV Climbing Daily. Why Does Climbing Have A Combined Format For The Olympics? Climbing Daily Ep.\ 1252. EpicTV Climbing Daily, September 19, 2018. \url{https:// www.youtube.com/watch?v=BdaBFobeQ7g}


\bibitem{Heywood} 
I.\ Heywood. (2006). Climbing monsters: excess and restraint in contemporary rock climbing. \emph{Leisure Studies} {\bf 25(4)} (2006),  455--467.


\bibitem{rules2018}
International Federation of Sport Climbing. \emph{IFSC Rules Modifications 2018}, April 2018.

\bibitem{rules2021}
International Federation of Sport Climbing. \emph{IFSC Rules Modifications 2021}, April 2021.

\bibitem{ioc}
International Olympic Committee. IOC approves five new sports for Olympic Games Tokyo 2020. International Olympic Committee, August 3, 2016. \url{https://olympics.com/ioc/news/ioc-approves-five-new-sports-for-olympic-games-tokyo-2020}

\bibitem{Kiewa}
J.\ Kiewa.   Traditional climbing: metaphor of resistance or metanarrative of oppression? \emph{Leisure Studies} {\bf 21(2)} (2002), 145--161.

\bibitem{Paris}
Sport Climbing. Paris 2024. \url{https://www.paris2024.org/en/sport/sport-climbing/}

\bibitem{plastic}
Plastic Weekly. Climbing Is Olympic, Camera Work Is Hard, Boulders Flop, Janja Ain't The GOAT Just Yet --- The Debrief. \emph{Plastic Weekly}, August 10, 2021. \url{https://www.youtube.com/watch?v=qxX0o8bgrSQ}


\bibitem{Samet}
M.\ Samet.  Upset In Sport Climbing Men's Speed Finals: Results Here. Climbing. August 5, 2021.\url{https://www.climbing.com/competition/olympics/upset-in-sport-climbing-mens-speed-finals-results-here/}

\bibitem{TW19} H.\ Thorpe, and B.\ Wheaton. The Olympic Games, Agenda 2020 and action sports: the promise, politics and performance of organisational change. \emph{International Journal of Sport Policy and Politics} {\bf11(3)} (2019), 465--483.

\bibitem{Putnam} William Lowell Putnam Mathematical Competition.
\url{https://en.wikipedia.org/wiki/William_Lowell_Putnam_Mathematical_Competition}.

\bibitem{Walker}
N.\ Walker. New Speed World Record at Salt Lake City World Cup. Gripped. May 29, 2021. \url{https://gripped.com/indoor-climbing/new-speed-world-record-at-salt-lake-city-world-cup/}

\end{thebibliography}
\end{document}